\documentclass[12pt]{amsart}
\usepackage{amscd}

\newtheorem{theorem}{Theorem}

\theoremstyle{definition}

\newtheorem{conjecture}{Conjecture}

\theoremstyle{remark}
\newtheorem*{remark}{Remark}

\usepackage{amsfonts}
\newcommand{\field}[1]{\mathbb{#1}}
\newcommand{\Q}{\field{Q}}
\newcommand{\R}{\field{R}}
\newcommand{\Z}{\field{Z}}

\newcommand{\bs}{\backslash}
\newcommand{\ra}{\rightarrow}

\begin{document}

\title[Discrete Components of some Complementary Series
Representations]{Discrete Components of Some Complementary Series
  Representations}

\author{B.Speh  and T. N. Venkataramana}

\email{venky@math.tifr.res.in, speh@math.cornell.edu}

\subjclass{Primary 11F75; Secondary  22E40, 22E41\\ B.Speh, Department
of  Mathematics,  310  Malott  Hall, Cornell  University,  Ithaca,  NY
14853-4201, U.S.A \\ T.  N. Venkataramana, School of Mathematics, Tata
Institute of Fundamental Research, Homi Bhabha Road, Bombay - 400 005,
INDIA.  venky@math.tifr.res.in}

\date{}

\begin{abstract}  We show  that the  restriction of  the complementary
reries  representations of  $SO(n,1)$ to  $SO(m,1)$ ($m<  n$) contains
complementary series representations of $SO(m,1)$ discretely, provided
that the continuous parameter is sufficiently close to the first point
of reducibility and the representation of $M$- the compact part of the
Levi- is a sufficiently small fundamental representation.  \\

We prove,  as a consequence, that the  cohomological representation of
degree  $i$ of  the group  $SO(n,1)$ contains  discretely,  for $i\leq
m/2$, the  cohomological representation of degree $i$  of the subgroup
$SO(m,1)$ if $i\leq m/2$. \\

As  a global  application,  we show  that  if $G/\Q$  is a  semisimple
algebraic group  such that $G(\R)=SO(n,1)$ up to  compact factors, and
if   we  assume  that   for  all   $n$,  the   tempered  cohomological
representations  are not limits  of complementary  series {\it  in the
automorphic  dual  of  $SO(n,1)$},  then  for  all  $n$,  non-tempered
cohomological representations are isolated  in the automorphic dual of
$G$.   This  reduces conjectures  of  Bergeron  to  the case  of  {\bf
tempered} cohomological representations.
\end{abstract}

\maketitle

\begin{flushright} \end{flushright}

\newpage
\section{Introduction}

A  famous  Theorem  of  Selberg  \cite{Sel}  says  that  the  non-zero
eigenvalues of the  Laplacian acting on functions on  quotients of the
upper half  plane $\mathfrak{h}$ by {\it congruence  subgroups} of the
integral modular group, are bounded  away from zero, as the congruence
subgroup  varies.   In  fact,   Selberg  proved  that  every  non-zero
eigenvalue  $\lambda$ of  the Laplacian  on functions  on  $\Gamma \bs
{\mathfrak h}$, $\Gamma \subset SL_2(\Z)$, satisfies the inequality
\[\lambda \geq \frac{3}{16}.\]

A  Theorem  of  Clozel  \cite{Clo}  generalises  this  result  to  any
congreunce quotient of any symmetric space of non-compact type: if $G$
is  a linear  semi-simple  group defined  over  $\Q$, $\Gamma  \subset
G(\Z)$ a  congruence subgroup and  $X=G(\R)/K$ the symmetric  space of
$G$, then nonzero eigenvalues $\lambda  $ of the Laplacian {\it acting
on the space of functions on} $\Gamma \bs X$ satisfy:
\[\lambda  \geq \epsilon,\]  where $\epsilon  >  0$ is  a number  {\bf
independent} of the congruence subgroup $\Gamma$. \\

Analogous  questions on  Laplacians  acting on  differential forms  of
higher degree  (functions may be  thought of as differential  forms of
degree  zero) have geometric  implications for  the cohomology  of the
locally symmetric space. Concerning  the eigenvalues of the Laplacian,
Bergeron (\cite{Ber}) has made the following conjecture:

\begin{conjecture}  (Bergeron) \label{bergeron}  Let $X$  be  the real
hyperbolic  $n$-space  and   $\Gamma  \subset  SO(n,1)$  a  congruence
arithmetic  subgroup.  Then  non-zero  eigenvalues $\lambda  $ of  the
Laplacian  acting  on  the  space   $\Omega  ^i  (\Gamma  \bs  X)$  of
differential forms of degree $i$ satisfy:
\[\lambda \geq \epsilon,\] for some $\epsilon >0$ {\bf independent} of
the congruence  subgroup $\Gamma$, provided $i$ is  strictly less than
the ``middle dimension (i.e. $i<[n/2]$).\\

If $n$  is even, the above  conclusion is conjectured to  hold even if
$i=[n/2]=n/2$. (When  $n$ is odd,  there is a slightly  more technical
statement which we omit).
\end{conjecture}

In  this paper, we  show, for  example, that  if the  above conjecture
holds true in the {\it middle  degree} for all even integers $n$, then
the conjecture holds for  differential forms of arbitrary degrees (See
Theorem  \ref{laplacian}).  For  odd  $n$,  there is,  again,  a  more
technical statement (see Theorem \ref{laplacian}). \\

The main Theorem  of the present paper is  Theorem \ref{mainth} on the
occurrence  of discrete components  in the  restriction of  of certain
complementary  series representations of  $SO(n,1)$ to  $SO(m,1)$. The
statement on  Laplacians may be  deduced from the main  theorem, using
the Burger-Sarnak method (\cite{Bu-Sa}). \\

We now describe the main  theorem more precisely.  Let $i\leq [n/2]-1$
and $G=SO(n,1)$.   Let $P=MAN$ be the Langlands  decomposition of $G$,
$K\subset G$  a maximal compact  subgroup of $G$ containing  $M$. Let
${\mathfrak p}_M$ be the standard representation of $M$ and $\wedge ^i$
be its  $i$-th exterior power  representation.  Denote by  $\rho _P^2$
the  character of  $P$ acting  on the  top exterior  power of  the Lie
algebra of $N$. Consider the representation
\[\widehat {\pi _u(i)}= Ind_P^G (\wedge ^i \otimes \rho _P(a)^u)\] for
$0<u<1-\frac{2i}{n-1}$.   The  representation  $\widehat {\pi  _u(i)}$
denotes  the {\bf  completion} of  the  space $\pi  _u$ of  $K$-finite
vectors with respect  to the $G$-invariant metric on  $\pi _u(i)$, and
is called the  {\bf complementary series representation} corresponding
to the representation $\wedge ^i$ of $M$ and the parameter $u$. \\

Let $H=SO(n-1,1)$ be embedded in $G$  such that $P\cap H$ is a maximal
parabolic  subgroup of  $H$, $A \subset  H$, and  $K\cap H$  a maximal
compact  subgroup   of  $H$.   We  now  assume   that  $\frac{1}{n-1}<
u<1-\frac{2i}{n-1}$.      Let      $u'=\frac{(n-1)u-1}{n-2}$;     then
$0<u'<1-\frac{2i}{n-2}$.  Denote by $\wedge ^i _H$ the $i$-th exterior
power of the standard representation of $M\cap H\simeq O(n-2)$.

We  obtain  analogously  the  complementary series representation 
\[\widehat  {\sigma _{u'}(i)}=  Ind_{P\cap H}^H  (\wedge ^i  _H \otimes
\rho _{P\cap H}(a)^{u'}),\]  of $H$. The main theorem  of the paper is
the following. 

\begin{theorem}\label{mainth} If 
\[\frac{1}{n-1}<u<1-\frac{2i}{n-1},\]  then  the complementary  series
representation  $\widehat{\sigma _{u'}(i)}$  occurs discretely  in the
restriction of  the complementary series  representation $\widehat{\pi
_u(i)}$ of $G=SO(n,1)$ to the subgroup $H=SO(n-1,1)$:
\[\widehat{\sigma _{u'}(i)}\subset \widehat{\pi _u(i)}_{|SO(n-1,1)}.\]
\end{theorem}

\begin{remark} The  corresponding statement is false for  the space of
$K$-finite vectors  in both the  spaces; this inclusion holds  only at
the level of completions of the representations involved.
\end{remark}

Denote by $A_j(n)$ the  unique unitary cohomological representation of
$G$ which has cohomology (with trivial coefficients) in degree $j$. As
$u$  tends  to   the  limit  $1-\frac{2i}{n-1}$,  the  representations
$\pi_u(i)$ tend (in  the Fell topology on the  Unitary dual ${\widehat
G}$)  both  to  the   representation  $A_i=A_i(n)$  and  to  $A_{i+1}=
A_{i+1}(n)$.  Using  this, and the  proof of Theorem  \ref{mainth}, we
obtain

\begin{theorem}  \label{cohomologicalreps}   The  restriction  of  the
cohomological  representation $A_i(n)$  of $SO(n,1)$  to  the subgroup
$H=SO(n-1,1)$  contains discretely,  the  cohomological representation
$A_i(n-1)$ of $SO(n-1,1)$:
\[A_i(n-1)\subset A_i(n)_{|SO(n-1,1)}.\]
\end{theorem}

Supose now  that $G$ is  a semi-simple linear algebraic  group defined
over $\Q$ and $\Q$-simple, such that
\[G(\R)\simeq SO(n,1)\] up to compact  factors (if $n=7$, we assume in
addition  that  $G$  is  not  the  inner  form  of  some  trialitarian
$D_4$). Then  there exists a $\Q$-simple  $\Q$-subgroup $H_1\subset G$
such that
\[H_1(\R)\simeq SO(n-2,1),\] up to compact factors.

Denote by $\widehat {G} _{Aut}$ the ``automorphic dual'' of $SO(n,1)$,
in the sense of  Burger-Sarnak (\cite{Bu-Sa}). Suppose $A_i=A_i(n)$ is
a limit  of representations $\rho  _m$ in $\widehat {G}  _{Aut}$.  The
structure of the unitary dual  of $SO(n,1)$ shows that this means that
$\rho _m=\widehat {\pi _{u_m}(i)}$ for some sequence $u_m$ which tends
from the left, to $1  - \frac{2i}{n-1}$ (or to $1 - \frac{2i+2}{n-1}$;
we   will  concentrate   only  on   the  first   case,  for   ease  of
exposition). \\

Since  $\rho  _m={\widehat  {\pi}}_{u_m}(i)\in  \widehat{G}_{Aut}$,  a
result of Burger-Sarnak (\cite{Bu-Sa}) (in the reference \cite{Bu-Sa},
the restriction  $\rho _m| _H$ refers  to the closure of  the union of
all the irreducibles of $H$ which occur weakly in $\rho _m$ restricted
to $H$, in  the Fell topology on the  unitary representations of $H$),
implies that
\[\rho _m |_H \subset \widehat{H}_{Aut}.\]

Applying Theorem \ref{mainth} twice, we get
\[\sigma_m={\widehat    \sigma}_{u_m''}(i)\in    \rho   _m |_H \subset 
{\widehat H}_{Aut}\]  Taking limits as  $m$ tends to infinity,  we get
$A_i(n-2)$   as    a   limit   of    representations   $\sigma_m$   in
$\widehat{H}_{Aut}$.    Therefore,the   isolation   of   $A_i(n)$   in
$\widehat{G}_{Aut}$ is reduced to  that for $SO(n-2,1),\cdots.$ We can
finally assume that $A_i(m)$ is a tempered representation of $SO(m,1)$
where $m=2i$ or $2i+1$. \\ This proves the following Theorem.

\begin{theorem}\label{laplacian}  If for all  $m$, the  {\bf tempered}
cohomological representation  $A_i(m)$ (i.e. when $i=[m/2]$)  is not a
limit of  complementary series in  the automorphic dual  of $SO(m,1)$,
then  for  all integers  $n$,  and  for  $i<[n/2]$, the  cohomological
representation  $A_i(n)$  is  isolated  in  the  automorphic  dual  of
$SO(n,1)$.
\end{theorem}

We  now sketch  a  proof of  the  Main Theorem.  Detailed proofs  will
apprear elsewhere.\\

The  proof of  Theorem  \ref{mainth} is  somewhat roundabout  and
proceeds as follows. 

(1)  We   first  prove  Theorem  \ref{mainth}  when   $i=0$;  that  is
$\widehat{\pi   }_u    =\widehat{\pi   }_u(i)$   is    an   unramified
representation and $\frac{1}{n-1}< u< 1$. In this case, we get a model
of the  representation $\widehat{\pi}_u$ by  restricting the functions
(sections of a line bundle) on  $G/P$ to the big Bruhat cell $Nw\simeq
\R  ^{n-1}$ and taking  their Fourier  transforms.  The  $G$- invariant
metric is  particularly easy  to work with  on the  Fourier
transforms;  it is  then easy  to  see that
$\widehat{\sigma}_{u'}$  embeds   isometrically  in
$\widehat{\pi}_u$. \\ 
(2) When this is interpreted
in the space $\pi_u$ of $K$-finite vectors, we have the isomorphism of
the intertwining map $I_G(u)$:
\[I_G  (u):  \pi_u  \simeq  \pi_{-u}=Ind_P ^G  (\rho  _P^{-u}).\]  The
restriction of the sections on  $G/P$ of the line bundle corresponding
to  $\pi   _{-u}$  to  the   subspace  $H/H\cap  P$  is   exactly  the
representation $\sigma_{-u'}$, which is isomorphic, via the inverse of
the intertwining  map, namely  $I_H(u')^{-1}$, to $\sigma  _{u'}$. The
existence    of     the    embedding     $\widehat{\sigma}_{u'}\subset
\widehat{\pi}_u$  (and  the multiplicity  one  for representations  of
$SO(n-1,1)$ occurring in $SO(n,1)$),  imply that the restriction map of
sections
\[res: \pi_{-u}\ra \sigma_{-u'},\]
is continuous for the metric on $\pi_{-u}\simeq \pi_u$ and on
$\sigma_{-u'}\simeq \sigma_{u'}$. \\

(2)  The  $K$-irreducible  representations  occurring in  $\pi_u$  are
parametrised  by  non-negative  integers  $m$, each  ireducible  $V_m$
occurring with multiplicity  one ($V_m$ is isomorphic to  the space of
homogeneous  harmonic   polynomials  of  degree  $m$   on  the  sphere
$G/P\simeq K/M$): $\pi_u=\oplus _{m\geq 0} V_m$.\\

Similarly,  we write  $\sigma_{u'}$ as  a direct  sum  of irreducibles
$W_l$  of   $K\cap  H$  (indexed  again   by  non-negative  integers):
$\sigma_{u'}=\oplus   _{l\geq  0}  W_l$.   Denote  by   $V_{m,l}$  the
isotypical component of $W_l$ in the restriction of $V_m$ to $K\cap H$
($V_m$ restricted to $K\cap H$  is in fact multiplicity free). We have
the  restriction  map  $\widehat{r}:  \pi= Ind_M^K(triv)  \ra  \sigma=
Ind_{M\cap H} ^{K\cap H}(Triv)$. This maps $V_{m,l}$ into $W_l$.  Set
\[C(m,l,0)=\frac{\mid\mid  \widehat{r}(f)   \mid\mid  ^2  _{K\cap  H}}
{\mid\mid f \mid\mid ^2 _{K}},\]

We show that the continuity of the map
$res: \pi_{-u}\ra \sigma_{-u'}$ is equivalent to the statement that
the series 
\begin{equation} \label{unramifiedconvergence} 
\sum _{m\geq l} C(m,l,0) \frac{l^{(n-2)u'}}{m^{(n-1)u}}< A,
\end{equation}
where $A$ is a constant independent of the integer $l$. 

The equivalence is proved by calculating the value of the intertwining
operator  on  each $K$-type  $V_m$  (the  operator  acts by  a  scalar
$\lambda_m(u)$ as a  ratio of values of the  classical Gamma function)
and obtaining asymptotics as $m$ tends to infinity by using Stirling's
Approximation to the classical Gamma function.

Since,  by the  result quoted  above,  the restriction  map is  indeed
continuous, the estimate  of equation (\ref{unramifiedconvergence}) is
indeed true. \\

(3) We now describe the  proof of Theorem \ref{mainth} in the ramified
case. We have analogously the restriction maps
\[r_u(i):\pi_{-u}(i)\ra       \sigma_{-u'}(i),\]      where      again
$u'=\frac{(n-1)u-1}{n-2}$.  By  the  multiplicity  one  statement  for
irreducible representations of $SO(n-1,1)$  which are quotients of the
representation $\pi_{-u}(i)$, it  follows that Theorem \ref{mainth} is
equivalent to the continuity of the restriction map $r_u(i)$. \\

As  a $K$-module,  again $\pi_u(i)$  is a  direct sum  of irreducibles
parametrised by  non-negative integers and  we can define  the numbers
$C(m,l,i)$ analogous  to the definition of $C(m,l,0)$  above. We prove
analogously  that  the  continuity  of  the  restriction  $r_u(i)$  is
equivalent to the estimate
\begin{equation} \label{ramifiedconvergence} 
\sum _{m\geq l} C(m,l,i)\frac{l^{(n-2)u'}}{m^{(n-1)u}}< A,
\end{equation}
where $A$ is a constant independent of the $K\cap H$ ``type'' $l$. The
proof of this equivalence uses some precise estimates of the
intertwining operators on the analogues of the spaces $V_m$ in
$\pi (i)\simeq  Ind_M^K(\wedge^i)$ (and on the analogues of $W_l$ in
$\sigma (i)$ which is defined similarly to $\pi (i)$ , but for the
group $H$).  \\ 

(4) We prove that the estimate in equation (\ref{ramifiedconvergence})
is indeed true, by analysing the $K$-types occurring in $\pi (i)$, and
showing as a consequence that
\[C(m,l,i)\leq  \gamma  C(m,l,0),\]  for  some constant  $\gamma $ 
independent  of
$m,l$.  Since equation  \ref{unramifiedconvergence}) is  proved  to be
true,  it  follows that  equation  (\ref{ramifiedconvergence}) also  is
true.  This proves Theorem \ref{mainth} for arbitrary $i$.\\

\noindent{\bf {Acknowledgments}.} \\ 

A considerable part of this work was done when T.N.V. was visiting the
Department  of Mathematics,  Cornell University  in 2007-2008  and the
Institute for Advanced Study, Princeton  in 2008. He thanks both these
Institutes  for their  hospitality; he  also thanks  Peter  Sarnak for
valuable conversations concerning the  material of the paper.  He also
gratefully acknowledges the support of the J.C Bose fellowship for the
period 2008-2013. \\

B. Speh was partially supported by NSF grant DMS 0070561. She would
alos like to thank the Tata Institute for its hospitality in the
winter of 2009 during a crucial stage of this work.

\end{document}